\documentclass{amsart}

\usepackage{amssymb,amsfonts, latexsym,amsmath,hhline,array,longtable}
\usepackage{pdfsync,color, comment,colortbl, mathrsfs,stmaryrd,cite,graphicx}
\usepackage[matrix,arrow,curve]{xy}

\usepackage{ifpdf}
\ifpdf \usepackage[colorlinks=true, citecolor=blue, linkcolor=blue, urlcolor=blue]{hyperref} \fi

\newcommand{\cal}{\mathcal}

\newtheorem{formula}{}[section]
\newtheorem{definition}[formula]{Definition}
\newtheorem{corollary}[formula]{Corollary}
\newtheorem{remark}[formula]{Remark}
\newtheorem{lemma}[formula]{Lemma}
\newtheorem{theorem}[formula]{Theorem}

\def\thrm{\begin{theorem}}
\def\thrml#1{\begin{theorem}\label{#1}}
\def\ethrm{\end{theorem}}
\def\rmrk{\begin{remark}}
\def\rmrkl#1{\begin{remark}\label{#1}}
\def\ermrk{\end{remark}}
\def\dfntn{\begin{definition}}
\def\dfntnl#1{\begin{definition}\label{#1}}
\def\edfntn{\end{definition}}
\def\nmrt{\begin{enumerate}}
\def\enmrt{\end{enumerate}}
\def\tm#1{\item[{\rm (#1)}]}
\def\qtnl#1{\begin{equation}\label{#1}}
\def\eqtn{\end{equation}}
\def\lmm{\begin{lemma}}
\def\lmml#1{\begin{lemma}\label{#1}}
\def\elmm{\end{lemma}}
\def\crllr{\begin{corollary}}
\def\crllrl#1{\begin{corollary}\label{#1}}
\def\ecrllr{\end{corollary}}
\def\css{\begin{cases}}
\def\ecss{\end{cases}}
						
\def\proof{\noindent{\bf Proof}.\ }

\def\cE{{\cal E}}

\def\cX{{\cal X}}

\def\fT{{\frak T}}
\def\fX{{\frak X}}

\DeclareMathOperator{\aut}{Aut}

\DeclareMathOperator{\AGL}{AGL}

\DeclareMathOperator{\dimwl}{dim_{\scriptscriptstyle WL}}

\DeclareMathOperator{\GL}{GL}

\DeclareMathOperator{\id}{id}

\DeclareMathOperator{\Inn}{Inn}

\DeclareMathOperator{\inv}{Inv}

\DeclareMathOperator{\syl}{Syl}
\DeclareMathOperator{\sym}{Sym}

\def\eprf{\hfill$\square$}

\def\phmb#1{{\phantom{x}\hspace{-2mm}^{#1}}}
\def\qaq{\quad\text{and}\quad}

\def\wt{\widetilde}

\begin{document}
	
\title{On pseudofrobenius imprimitive association schemes}
\author{Ilia Ponomarenko}
\address{St. Petersburg Department of V.A. Steklov Institute of Mathematics, St. Petersburg, Russia}
\email{inp@pdmi.ras.ru}
\author{Grigory Ryabov}
\address{Sobolev Institute of Mathematics, Novosibirsk, Russia}
\address{Novosibirsk State Technical University, Novosibirsk, Russia}
\address{St. Petersburg Department of V.A. Steklov Institute of Mathematics, St. Petersburg, Russia }
\address{Leonard Euler International Mathematical Institute in Saint Petersburg, St. Petersburg, Russia}
\email{gric2ryabov@gmail.com}
\thanks{The second author is supported by Leonard Euler International Mathematical Institute in Saint Petersburg under agreement No.~075-15-2019-1620 with the Ministry of Science and Higher Education of the Russian Federation}
\date{}
	
\begin{abstract}
An (association) scheme is said to be \emph{Frobenius} if it is the scheme of a Frobenius group. A scheme which has the same tensor of intersection numbers as some Frobenius scheme is said to be \emph{pseudofrobenius}. We establish a necessary and sufficient condition for an imprimitive pseudofrobenius scheme to be Frobenius. We also prove strong necessary conditions for existence of an imprimitive pseudofrobenius scheme which is not Frobenius. As a byproduct, we obtain a sufficient condition for an imprimitive Frobenius group $G$ with abelian kernel to be determined up to isomorphism only by the character table of~$G$. Finally, we prove that the Weisfeiler-Leman dimension of a circulant graph with~$n$ vertices and Frobenius automorphism group is equal to~$2$ unless $n\in \{p,p^2,p^3,pq,p^2q\}$, where $p$ and $q$ are distinct primes. 
\end{abstract}
	
\maketitle

\section{Introduction}

Frobenius groups play an important role in the group theory, geometry, and combinatorics. Every Frobenius group has the standard permutation representation in which the Frobenius kernel acts regularly and a one-point stabilizer coincides with  Frobenius complement. The corresponding permutation groups often act as automorphism groups of geometrical (skew-affine spaces~\cite{Andre}) or combinatorial (association schemes) objects. For example, every imprimitive Frobenius group is the automorphism group of a Frobenius association scheme~\cite{MP}; for the basics of (association) schemes, see Section~$3$. A natural question arises: is it possible to describe Frobenius groups in the framework of ``group theory without groups''~\cite{BI}, i.e., in a pure combinatorial way. The main goal of this paper is to give a partial answer to this question for imprimitive Frobenius groups.\medskip

In attempt to answer this question, one faces with the following problem: there are schemes which are combinatorially similar (i.e., have  the same tensor of intersection numbers) but only one of them is Frobenius. A scheme which has the same tensor of intersection numbers as the scheme of some Frobenius group~$G$ is said to be \emph{$G$-pseudofrobenius}, or just \emph{pseudofrobenius} if $G$ is unessential. Every pseudofrobenius scheme belongs to the well-known class of pseudocyclic schemes studied in, e.g., in~\cite{BCN,CHPV,MP}.\medskip


Asymptotically every pseudofrobenius scheme is Frobenius~\cite[Theorem~1.1]{MP}. On the other hand, there are exponentially many (imprimitive) pseudofrobenius schemes which are not Frobenius (see below). All of these schemes are not schurian which means that they are not the schemes of permutation groups. In contrast to this, our first main result shows that if a pseudofrobenius scheme is ``locally schurian'', then it is Frobenius. The property to be ``locally schurian'' is expressed in the language of the $4$-condition introduced by Hestenes and Higman for strongly regular graphs~\cite{HH} and later generalized to association schemes (see, e.g.,~\cite{MP}).

\thrml{4condition}
An imprimitive pseudofrobenius scheme is Frobenius if and only if it satisfies the $4$-condition. 
\ethrm

It should be mentioned that a special case of Theorem~\ref{4condition} for some class of imprimitive Frobenius schemes was proved in~\cite[Theorem~4.1]{MP}. Further, Theorem~\ref{4condition} substantially improves~\cite[Theorem~5.11]{EP} which states, in fact, that an imprimitive pseudofrobenius scheme is Frobenius whenever it satisfies a condition stronger than $6$-condition. \medskip

The second result of the paper establishes a strong necessary condition for the existence of pseudofrobenius schemes which are not Frobenius; in what follows they are called \emph{proper} pseudofrobenius schemes. To state this result, we need some preparation. Let $G$ be a Frobenius group, and $H$ and $K$ the kernel and complement of $G$, respectively. Denote by~$\pi(H)$ the set of all prime divisors of $|H|$ and by~$d(G)$ the length of a maximal chain of $K$-invariant subgroups of $H$. Note that $d(G)$ does not depend on the choice of such chain by~\cite[Theorem~6]{Bourbaki}. A \emph{principal section} of~$G$ is defined to be the primitive permutation group induced by $G$ on the set of all right $L$-coset inside $U$, where $L$ and $U$ are successive  members of the longest chain of $K$-invariant subgroups of $H$ (see Section~2 for the exact definitions).

\thrml{2}
Let $G$ be an imprimitive Frobenius group with kernel $H$, $\pi=\pi(H)$, and $d=d(G)$. Suppose that there exists a proper $G$-pseudofrobenius scheme. Then $(|\pi|,d)\in \{1,2\}\times \{2,3\}$. Moreover, if $d=3$, then

\nmrt

\tm{1} $|H|=p^{3a}$ for a prime $p$ and an integer~$a\geq 1$, and every principal section of~$G$ is $2$-transitive of degree~$p^a$, or

\tm{2}  $|H|=p^{2a}q^b$ for distinct primes $p$ and $q$ and integers~$a,b\geq 1$, and every principal section of $G$ is either $2$-transitive of degree $p^a$ or of rank~$3$ and degree~$q^b$; in particular, $q^b=2p^a-1$.

\enmrt 

\ethrm

Theorem~\ref{2} is a consequence of Theorem~\ref{separability} stating that a pseudofrobenius scheme is determined up to isomorphism by the tensor of its intersection numbers unless the conclusion of Theorem~\ref{2} holds. The last result weakens the condition of~\cite[Theorem~1.3]{MP} and strengthens the conclusion of~\cite[Theorem~1.1]{Ry}.\medskip

Let us give a well-known example of a proper pseudofrobenius scheme. Let $q$ be a prime power, $H$ the translation group of a $2$-dimensional vector space over the finite field $\mathbb{F}_q$, and $K$ the center of $\GL(2,q)$. Then the semidirect product $G=HK\leq \AGL(2,q)$ is an imprimitive Frobenius group with kernel~$H$, complement $K$, and $(|\pi|,d)=(1,2)$. It is known (see, e.g., \cite[Section~2.5.2]{CP}), that there are proper $G$-pseudofrobenius schemes. We do not know whether there exists a proper $G$-pseudofrobenius scheme for an imprimitive Frobenius group~$G$ with $(|\pi|,d)\in \{(1,3),(2,2),(2,3)\}$. Moreover, if $d=3$ and $|H|\leq 3^4\,17$, then a proper $G$-pseudofrobenius scheme does not exist (this was verified by computer calculations).\medskip

By Thompson's theorem, the kernel of every Frobenius group is nilpotent. The corollary below establishes the properties of the kernel of a Frobenius group from Theorem~\ref{2}.

\crllrl{nilpotent}
In the conditions of Theorem~\ref{2}, the nilpotency class~$c$ of $H$ is at most~$3$. Moreover, suppose that either $H$ is not a $2$-group or $d\neq 3$. Then
\nmrt

\tm{1} $c\leq 2$ and $c=2$ only if $(|\pi|,d)=(1,2)$, or $(|\pi|,d)=(2,3)$ and $2\in \pi$; 

\tm{2} $H$ is abelian if $|H|$ is odd and $(|\pi|,d)\neq (1,2)$.

\enmrt
\ecrllr

Corollary~\ref{nilpotent} implies that if $(|\pi|,d)=(2,2)$, then the kernel of the group $G$ is abelian, and there are infinitely many such groups. In the remaining cases, namely, if $(|\pi|,d)\in \{(1,2),(1,3),(2,3)\}$, there are also infinitely many examples with cyclic kernels. Concerning nonabelian case, one can take the normalizer of a Sylow $2$-subgroup of a Suzuki group; in this case the kernel of $G$ is of nilpotency class~$2$ and order~$2^{2m}$, where $m$ is odd, and the complement is cyclic of order~$2^m-1$. This gives an example with $(|\pi|,d)=(1,2)$ (see~\cite{Suz}). Finally, in the cases $(|\pi|,d)=(1,3)$ or~$(2,3)$, the corresponding examples can be constructed as subdirect products of the above example and a Frobenius group with elementary abelian kernel of order~$2^{m}$ or cyclic kernel of prime order~$2^{m+1}-1$ (the Mersenne prime), respectively.\medskip

With the help of  technique used  to prove  Theorem~\ref{2}, we derive in Section~8 the following statement showing that there are no Brauer pairs among imprimitive Frobenius groups with abelian kernel, for which the conclusion of Theorem~\ref{2} does not hold.

\thrml{3}
Let $G_1$ and $G_2$ be imprimitive Frobenius groups with abelian kernels. Suppose that the conclusion of Theorem~\ref{2} does not hold for $G_1$. Then $G_1$ and $G_2$ have the same character tables if and only if $G_1$ and $G_2$ are isomorphic.
\ethrm

The technique developed in the present paper enables us to find the Weisfeiler-Leman dimension of many Frobenius circulants. Recall that a \emph{circulant} is a graph~$\fX$ isomorphic to a Cayley graph over cyclic group; $\fX$ is said to be \emph{Frobenius} if the group $\aut(\fX)$ is Frobenius. Among the examples of such graphs, there are connected circulants with $p^m$ vertices, where $p$ is an odd prime and $m\geq 1$, and connection set of cardinality at most~$p-1$. The \emph{Weisfeiler-Leman dimension} $\dimwl(\fX)$ of a graph~$\fX$ can be thought as the minimum number of variables in a formula of a natural fragment of first-order logic with counting quantifiers, which is valid only for graphs isomorphic to~$\fX$~\cite{GN}. Using~\cite[Theorem~4.4.7]{CP}, it is not hard to check that $\dimwl(\fX)\leq 3$ for every Frobenius circulant $\fX$.

\thrml{circulant}
Let $\fX$ be a Frobenius circulant on $n$ vertices. Then $\dimwl(\fX)=2$ unless $n\in \{p,p^2,p^3,pq,p^2q\}$, where $p$ and $q$ are distinct primes.
\ethrm

Observe that even in the exceptional cases of Theorem~\ref{circulant}, there are infinitely many Frobenius circulants $\fX$ with $\dimwl(\fX)=2$. On the other hand, there are Frobenius circulants~$\cX$ with prime number of vertices and $\dimwl(\fX)=3$ (see, e.g.,~\cite[Section~4.5]{CP}). However, we do not know such circulants with composite number of vertices.

The authors are grateful to the participants of the seminar ``Discrete Algebraic Structures and Isomorphism Problem'' for the fruitful discussions on the subject matters and the valuable comments.

\section{Frobenius groups}

Throughout this section, $G$ is a Frobenius group and $H$ and $K$ are the Frobenius kernel and complement of $G$, respectively. The group~$H$ is nilpotent by the Thompson theorem (see, e.g.,~\cite[Theorem~6.24]{Is}). In particular, every Sylow subgroup of~$H$ is characteristic and hence $K$-invariant. Denote by $\mathcal{H}(G)$ the set of all $K$-invariant subgroups of $H$. Clearly, $\mathcal{H}(G)$ forms a sublattice of the subgroup lattice of $H$.\medskip  

Throughout the paper, we consider $G$ as a permutation group on~$H$, where $H$ acts on itself by right multiplications whereas $K$ acts by conjugations. In particular,~$H$ is a regular subgroup and~$K$ is a one-point stabilizer of~$G$. Thus, the following lemma is an easy consequence of~\cite[Chapter~IV]{Wi}.

\lmml{invariant}
A partition of $H$ is imprimitivity system of $G$ if and only if it is the partition into the right $L$-cosets for some $L\in \mathcal{H}(G)$.
\elmm



Let us recall the concept of section of a transitive group~\cite[Section~9]{EP2}. To this end, we define a \emph{factor block} $S$ of $G$ as the set of all blocks of an imprimitivity system of $G$, contained in a block $X$ of a smaller imprimitivity system.\footnote{An imprimitivity system $I$ is smaller than an imprimitivity system $I'$ if each block of $I'$ is contained in some block of $I$.} Clearly, the setwise stabilizer of $X$ in $G$ acts naturally on $S$. The permutation group $G^S\leq \sym(S)$ induced by this action is called the \emph{section} of a transitive group $G$. Note that if $S$ and $T$ are factor blocks corresponding to the same imprimitivity system, then the groups $G^S$ and $G^T$ are permutation isomorphic.\medskip

Now let $G$ be a Frobenius group as at the beginning of the section. Then by Lemma~\ref{invariant}, every group $U\in \mathcal{H}(G)$ is a block of $G$. In what follows, we deal with factor blocks of the form $S=U/L$ (the set of all right $L$-cosets inside $U$), where $L\in \mathcal{H}(G)$; if $L=\{1\}$, then we identify $S$ with $U$. We say that a section $G^S$ is \emph{principal} if there exists a chain of $K$-invariant subgroups 
$$\{1\}=H_0<H_1<\ldots <H_d=H$$ 
with $d=d(G)$ and such that $S=H_{i+1}/H_i$ for some $i\in\{0,\ldots,d-1\}$. From the definition of $d(G)$ and Lemma~\ref{invariant}, it easily follows that every principal section of~$G$ is primitive.

\lmml{class}
In the above notations,
\nmrt

\tm{1} $c(H)\leq d(G)$;

\tm{2} if $U\in \mathcal{H}(G)$ and $U\neq \{1\}$, then $G^U$ is Frobenius, and also $d(G^U)\leq d(G)-1$ if $U\neq G$;

\tm{3} every primitive section of $G$ has a prime power degree;

\tm{4} if $|K|$ is even, then $H$ is abelian.

\enmrt

\elmm

\proof
Every group of the upper central series of $H$ being characteristic in~$H$, is $K$-invariant. Consequently, $c(H)\leq d(G)$ and statement~$(1)$ holds. Statement~$(2)$ is obvious. Let $S=U/L$ be a factor block such that $G^S$ is primitive. Since the group $G^U\leq \sym(U)$ is Frobenius, its kernel is a normal nilpotent subgroup. Therefore, the group $G^S$ has a normal nilpotent subgroup. By the O'Nan-Scott theorem, this implies that the socle of $G^S$ is elementary abelian and regular and hence $G^S$ has a prime power degree. Thus, statement~$(3)$ holds. Statement~$(4)$ follows from~\cite[Theorem~6.3]{Is}.
\eprf

\section{Association schemes}
In our presentation of (association) schemes, we follow monograph~\cite{CP}; the proofs of all facts not explained below can be found there.

\subsection{Definitions}

Let~$\Omega$ be a finite set and~$S$ a partition of~$\Omega \times \Omega$; the elements of~$S$ are treated as binary relations on~$\Omega$. For arbitrary~$\alpha \in \Omega$ and $s\in S$, we put $\alpha s=\{\beta\in \Omega:\ (\alpha,\beta)\in s\}$. The pair~$\mathcal{X}=(\Omega,S)$ is called an (\emph{association}) \emph{scheme} on~$\Omega$ if the following conditions are satisfied:
\nmrt
\tm{C1}  the diagonal relation $1_{\Omega}=\{(\alpha,\alpha):\ \alpha\in\Omega\}$ belongs to~$S$;
\tm{C2} for each $s\in S$, the relation $s^*=\{(\alpha,\beta): (\beta,\alpha)\in s\}$ belongs to~$S$;
\tm{C3}  given $r,s,t\in S$, the number $c_{rs}^t=|\alpha r\cap \beta s^{*}|$ does not depend on $(\alpha,\beta)\in t$. 
\enmrt

Any relation belonging to $S$ is called a \emph{basis relation} of~$\cX$; the basis relation containing a pair $(\alpha,\beta)$ is denoted by $r(\alpha,\beta)$. The numbers $|\Omega|$ and~$c_{rs}^t$ are called the \emph{degree} and \emph{intersection numbers} of~$\cX$, respectively. Given $r,s\in S$, the set $\{t\in S:~c_{rs}^t>0\}$ is denoted by $rs$.\medskip

Let $s\in S$. For every $\alpha\in \Omega$, the number $n_s=|\alpha s|$ called the \emph{valency} of $s$ is equal to~$c_{ss^*}^{1_{\Omega}}$ and hence does not depend on the choice of $\alpha$. It is well-known that
\qtnl{triangle}
n_tc_{rs}^{t^*}=n_rc_{st}^{r^*}=n_sc_{tr}^{s^*}
\eqtn
for all $r,s,t\in S$ (see, e.g.,~\cite[Eq.~2.1.14]{CP}).\medskip

An equivalence relation $e$ on the set $\Omega$ is called a \emph{parabolic} of the scheme~$\cX$ if~$e$ is the union of some basis relations of $\cX$. Clearly, $1_{\Omega}$ and $\Omega \times \Omega$ are parabolics of~$\cX$; we call them \emph{trivial}. A scheme is said to be \emph{imprimitive} if it has a nontrivial parabolic and \emph{primitive} otherwise. The set~$\cE(\cX)$ of all parabolics of $\cX$ forms a sublattice of the lattice of equivalence relations on $\Omega$. 

\subsection{Isomorphisms and schurity}

Let $\cX=(\Omega,S)$ and $\cX'=(\Omega',S')$ be schemes. A bijection $f:\Omega \rightarrow \Omega'$ is called a (\emph{combinatorial}) \emph{isomorphism} from $\mathcal{X}$ to $\mathcal{X}'$ if $S^f=S'$, where $S^f=\{s^f:~s\in S\}$ and $s^f=\{(\alpha^f,\beta^f):~(\alpha,\beta)\in s\}$. The set of all isomorphisms from $\cX$ to itself is a permutation group on $\Omega$ which has a normal subgroup 
$$\aut(\cX)=\{f\in \sym(\Omega):~s^f=s~\text{for all}~s\in S\},$$
called the \emph{automorphism group} of $\cX$.\medskip

Let $G$ be a transitive subgroup of $\sym(\Omega)$. Denote by $(\alpha,\beta)^G$ the orbit of the induced action of $G$ on $\Omega\times \Omega$, which contains  the pair $(\alpha,\beta)$. Then
$$\inv(G)=\inv(G,\Omega)=(\Omega, \{(\alpha,\beta)^G:~\alpha,\beta\in \Omega\})$$
is a scheme called the \emph{scheme of}~$G$. A scheme~$\cX$ is said to be \emph{schurian} if~$\cX$ is a scheme of some transitive permutation group.

\subsection{Algebraic isomorphisms and separability}

A bijection $\varphi:S\rightarrow S'$ is called an \emph{algebraic isomorphism} from $\cX$ to $\cX'$ if $c_{r^\varphi s^\varphi}^{t^\varphi}=c_{rs}^t$ for all $r,s,t\in S$. It is known that $|\Omega|=|\Omega'|$, $1^{\varphi}_{\Omega}=1^{}_{\Omega'}$, and $n_{s^{\varphi}}=n_s$ for every $s\in S$. The algebraic isomorphism $\varphi$ is extended to the unions of basis relations, namely, if $s_1,s_2\ldots\in S$, then we put $(s_1\cup s_2\cup\cdots)^{\varphi}=s_1^\varphi\cup s_2^\varphi\cup\cdots$. Note that if $e\in \mathcal{E}(\mathcal{X})$, then $e^{\varphi}\in \mathcal{E}(\mathcal{X}')$.\medskip

Every isomorphism $f$ from $\cX$ to $\cX'$ induces the algebraic isomorphism from $\cX$ to $\cX'$ which maps $s\in S$ to $s^f\in S'$. A scheme is said to be \emph{separable} if every algebraic isomorphism from it to another scheme is induced by an isomorphism.

\subsection{The $t$-condition} 

The following definition goes back to~\cite{FKM} (see also~\cite[Section~2.7]{MP}). Let $t\ge 2$ be an integer and $\fT(\cX)=\{T(\gamma):\ \gamma\in\Omega^t\}$, where  $T(\gamma)$  is a $t\times t$ array defined as follows:
\qtnl{241020z}
T(\gamma)_{ij}=r(\gamma_i,\gamma_j),\quad 1\le i,j\le t,
\eqtn
where $\gamma_i$ and $\gamma_j$ are the $i$th and $j$th coordinates of $\gamma$, respectively. We say that the scheme $\cX$ satisfies the {\it $t$-condition}  if for all $r\in S$ and $T\in \fT(\cX)$, the number
$$
c_T^r(\alpha,\beta)=|\{\gamma\in\Omega^t:\ T(\gamma)=T,\ (\gamma_1,\gamma_2)=(\alpha,\beta)\}|
$$ 
does not depend on $(\alpha,\beta)\in r$; we denote this number by $c_T^r$. Note that $\cX$ always satisfies the $3$-condition (and, of course, $2$-condition); in this case, the numbers~$c_T^r$ are just the intersection numbers of~$\cX$. A scheme is schurian if and only if it satisfies the $t$-condition for every $t\geq 2$~\cite[p.~5]{MP}.

\section{Pseudofrobenius schemes} 

Let $G$ be a Frobenius group. Following~\cite{MP}, the scheme $\inv(G)$ is said to be a \emph{Frobenius}. The parabolics of $\inv(G)$ are exactly $G$-invariant equivalence relations. The latter are in one-to-one correspondence with the imprimitivity systems of $G$: the classes of a $G$-invariant equivalence relation are the blocks of the corresponding imprimitivity system. By Lemma~\ref{invariant}, the imprimitivity systems of $G$ are the partitions into the right cosets by a subgroup from $\mathcal{H}(G)$. Thus, we have a bijection
\qtnl{bijection}
\chi: \mathcal{E}(\inv(G))\rightarrow \mathcal{H}(G)
\eqtn 
between the parabolics of $\inv(G)$ and groups belonging to~$\mathcal{H}(G)$.\medskip

A scheme $\mathcal{X}=(\Omega,S)$ is said to be \emph{pseudofrobenius} if it is algebraically isomorphic to the scheme $\inv(G)$, where $G$ is a (not necessary uniquely determined) Frobenius group; in this case the scheme $\cX$ is called \emph{$G$-pseudofrobenius}. Clearly, every Frobenius scheme is pseudofrobenius. The converse statement is not true in general, see Introduction. A pseudofrobenius scheme is said  to be \emph{proper} if it is not isomorphic to a Frobenius scheme. Examples of proper pseudofrobenius schemes can be found in~\cite[Sections~2.5.2,~4.5]{CP}.

\lmml{parabolics} 
Let $G$ be a Frobenius group, $\cX$ a $G$-pseudofrobenius scheme, and $\varphi$ an algebraic isomorphism from $\inv(G)$ to $\cX$. Then 
\nmrt
\tm{1} the degree of $\cX$ and valency of an irreflexive basis relation of $\cX$ are equal to the orders of the kernel and complement of $G$, respectively,

\tm{2} the mapping
$$\mathcal{H}(G)\rightarrow \mathcal{E}(\mathcal{X}),~L\mapsto e(L),$$
where $e(L)=L^{\chi^{-1}\varphi}$ with $\chi$ being bijection~\eqref{bijection}, is a lattice isomorphism such that $|L|=n_{e(L)}$. 

\enmrt 
\elmm

\proof
Statement~(1) is obvious. To prove statement~(2), we make use of some simple facts proved in monograph~\cite{CP}; the references to them are given without mentioning~\cite{CP}.

By~Exercise~1.4.16~(6), the bijection $\chi^{-1}$ coincides with the restriction of the bijection $\rho$ constructed in~Exercise~1.4.15 to $\mathcal{H}(G)$. By Exercise~1.4.15~(4),~(5), the mapping $\chi^{-1}:\mathcal{H}(G)\rightarrow \mathcal{E}(\inv(G))$ is a lattice isomorphism. On the other hand, the bijection $\mathcal{E}(\inv(G))\rightarrow \mathcal{E}(\cX)$, $e\mapsto e^{\varphi}$ is also a lattice isomorphism by~Exercise~2.7.30. Thus, so is the composition of $\chi^{-1}$ and $\varphi$. Moreover for every $L\in \mathcal{H}(G)$, we have
$$|L|=n_{L^{\chi^{-1}}}=n_{L^{\chi^{-1}\varphi}},$$
where the first equality follows from~formula~1.4.15 whereas the second one follows from~Corollary~2.3.20.
\eprf\medskip

By statement~(1) of Lemma~\ref{parabolics}, the pseudofrobenius scheme $\cX$ is \emph{equivalenced}, i.e., the number $n_s$ does not depend on irreflexive $s\in S$; it is called the \emph{valency} of~$\mathcal{X}$. Now as an easy consequence of Eq.~\eqref{triangle}, we obtain
\qtnl{trianglefrob}
c_{rs}^{t^*}=c_{st}^{r^*}=c_{tr}^{s^*},~r,s,t\in S\setminus \{1_{\Omega}\}.
\eqtn

\lmml{intersection}
Let $\mathcal{X}$ be an imprimitive pseudofrobenius scheme, $e$ a nontrivial parabolic of $\mathcal{X}$, and relations $r,s,t\in S$ irreflexive. Assume that 
$$t\in rs,\quad t\subseteq e,\quad r\cup s\not\subseteq e.$$ 
Then $c_{rs}^t=1$ and $c_{rs}^u=0$ for every $u\in S$ such that $u\subseteq e$ and $u\neq t$.
\elmm

\proof
We note that $c_{rs}^t\geq 1$ and hence $c_{st^*}^{r^*}\geq 1$ and $c_{t^*r}^{s^*}\geq 1$ (see Eq.~\eqref{trianglefrob}). So $r^*\in st^*$ and $s^*\in t^*r$. This implies that if one of $r$ and $s$ is contained in $e$, then, because of $t\subseteq e$, so is the other one. However, this contradicts the condition $r\cup s\not\subseteq e$. Therefore both of $r$ and $s$ lie outside $e$.  

Assume on the contrary that either $c_{rs}^t>1$ or there exists $u\in S$ such that $u\subseteq e$, $u\neq t$, and $c_{rs}^u>0$. Then $c_{st^*}^{r^*}>1$ or $c_{su^*}^{r^*}>0$ by Eq.~\eqref{trianglefrob}. In both cases for every $(\alpha,\beta) \in r^*$, there are distinct $\gamma,\delta\in \Omega$ such that 
$$(\alpha,\gamma),(\alpha,\delta)\in s~\text{and}~(\gamma,\beta),(\delta,\beta)\in e.$$ 
The second part of this formula implies that $\beta,\gamma,\delta$ belong to the same class $\Delta$ of the parabolic $e$. Therefore by the first part of this formula, $|\alpha s\cap \Delta|\geq 2$. Moreover, $\alpha\notin \Delta$ for otherwise $r=r(\beta,\alpha)\subset e$ which is impossible by the first paragraph. Now we arrive to a contradiction with~\cite[Eq.~(3.3.3)]{CP} stating that every point $\alpha$ of an equivalenced scheme has at most one $s$-neighbor in a class of $\Delta \not\ni \alpha$ of any nontrivial parabolic~$e$.
\eprf\medskip

 We complete the section by a statement which will be used in Section~7.

\lmml{divide}
Let $\cX$ be a pseudofrobenius scheme and $e_1\subseteq e_2$ parabolics of $\mathcal{X}$. Then the valency of $\mathcal{X}$ divides $\frac{n_{e_2}}{n_{e_1}}-1$.
\elmm

\proof
Denote by $G$ the Frobenius group such that $\cX$ is $G$-pseudofrobenius. Let $K$ be a complement of $G$, and let $L_1,L_2\in \mathcal{H}(G)$ be such that $e_1=e(L_1)$, $e_2=e(L_2)$, see Lemma~\ref{parabolics}(2). Then the set $L_2\setminus L_1$ is $K$-invariant. Therefore, the valency $k=|K|$ of $\cX$ divides the number 
$$|L_2|-|L_1|=|L_1|\,\left(\frac{|L_2|}{|L_1|}-1\right).$$ 
On the other hand, $|L_1|$ is coprime to~$k$, because the group $G^{L_1}$ is Frobenius or trivial by Lemma~\ref{class}(2). Thus, $k$ divides $\frac{|L_2|}{|L_1|}-1=\frac{n_{e_2}}{n_{e_1}}-1$.
\eprf

\section{Base triples and algebraic isomorphisms}

Throughout the section, $\cX=(\Omega,S)$ is an imprimitive pseudofrobenius scheme and $e$ is a nontrivial parabolic of $\cX$. 

\subsection{Base triples}
Let $\mu,\nu,\rho$ be points of~$\cX$.  We say that  $\tau=(\mu,\nu,\rho)$ is  a {\it base triple}  with respect to the parabolic~$e$ if $\mu\neq \nu$ and 
\qtnl{200920a}
(\mu,\nu)\in e\qaq (\mu,\rho)\not\in e.
\eqtn
Clearly, $\mu$, $\nu$, and $\rho$ are pairwise distinct. Denote by $S_\tau$ the set of all pairs $(x,y)\in S\times S$ in which $x$ is arbitrary and $y$ satisfies the following conditions:
\qtnl{200920w}
xy^*\ni  \css
r(\mu,\rho) &\text{if $x\subseteq e$,}\\
 r(\mu,\nu) &\text{if $x\nsubseteq e$.}\\
 \ecss
\eqtn
In view of Eq.~\eqref{200920a}, we have $y\not\subseteq e$. It should also be noted that $S_\tau$ depends only on the relations $e$, $r(\mu,\rho)$, and $r(\mu,\nu)$, and does not depend on the triple~$\tau$.\medskip

For each point $\alpha\in\Omega$, we define its coordinates $(x_\alpha,y_\alpha)\in S_\tau$ with respect to the base triple~$\tau$ as follows:
 \qtnl{200920w3}
 x_\alpha=r(\mu,\alpha)\qaq y_\alpha=
 \css
 r(\rho,\alpha) &\text{if $(\mu,\alpha)\in e$,}\\
 r(\nu,\alpha) &\text{if $(\mu,\alpha)\not\in e$}.\\
 \ecss
 \eqtn
 The following lemma shows that $\alpha$ is uniquely determined by its coordinates.

\lmml{170920c}
In the above notation, the mapping
$$
f_\tau:\Omega\to S_\tau,\ \alpha \mapsto (x_\alpha,y_\alpha).
$$
is a bijection.  Moreover, setting $r=r(\mu,\rho)$, $s=r(\mu,\nu)$, and $t=r(\nu,\rho)$, we have
\qtnl{130920a}
\mu^{f_\tau}=(1_{\Omega},r^*),\quad \nu^{f_\tau}=(s,t^*),\quad \rho^{f_\tau}=(r,t).
\eqtn
\elmm
\proof Eq.~\eqref{130920a} is straightforward consequence of Eq.~\eqref{200920w3}. Let $(x,y)\in S_\tau$. Let us prove that there exists $\alpha \in \Omega$ such that $(x_{\alpha},y_{\alpha})=(x,y)$. First assume that $x\subseteq e$. Since $y\nsubseteq e$ and $r\nsubseteq e$, we have $c_{y^*r^*}^{x^*}=1$ by Lemma~\ref{intersection}. So $c_{xy^*}^{r}=c_{y^*r^*}^{x^*}=1$ by Eq.~\eqref{trianglefrob}. Therefore,
$$
\mu x\cap \rho y=\{\alpha\}
$$  
for some $\alpha\in \Omega$. Thus, $x=r(\mu,\alpha)=x_{\alpha}$ and $y=r(\rho,\alpha)=y_{\alpha}$. Now let $x\nsubseteq e$. Since $y^*\nsubseteq e$ and $s\subseteq e$, we have $c_{xy^*}^s=1$ by Lemma~\ref{intersection}. Therefore,
$$
\mu x\cap \nu y=\{\alpha\}
$$  
for some $\alpha\in \Omega$. Thus, $x=r(\mu,\alpha)=x_{\alpha}$ and $y=r(\nu,\alpha)=y_{\alpha}$. This shows that the mapping $f_\tau$ is a surjective.\medskip

To prove that $f_\tau$ is injective, assume on the contrary that there are  distinct $\alpha,\beta\in\Omega$ such that $(x_\alpha,y_\alpha)=(x_\beta,y_\beta)=:(x,y)$. Then $\alpha,\beta\in\mu x$. If $x\subseteq e$, then $\alpha,\beta\in\rho y$ and hence  $c_{xy^*}^r=c_{y^*r^*}^{x^*}\ge 2$ which contradicts Lemma~\ref{intersection}. Finally, if $x \nsubseteq e$, then $\alpha,\beta\in\nu y$, and hence  $c_{xy^*}^s\ge 2$ which contradicts Lemma~\ref{intersection}. \eprf

\subsection{Invariance.} In what follows, we fix a base triple $\tau$ of the scheme $\cX$ with respect to the parabolic~$e$. Let $\cX'=(\Omega',S')$ be a scheme and $\varphi:S\to S'$, $s\to s'$ algebraic isomorphism. The following statement is straightforward.

\lmml{221020a}
In the above notation, $\cX'$ is a pseudofrobenius scheme and $e'=\varphi(e)$ a parabolic of~$\cX'$. Moreover, there exists a base triple $\tau'=(\mu',\nu',\rho')$ of the scheme~$\cX'$ with respect to the parabolic~$e'$, such that 
\qtnl{181020b}
r(\mu,\nu)'=r(\mu',\nu'),\quad r(\mu,\rho)'=r(\mu',\rho'),\quad r(\nu,\rho)'=r(\nu',\rho').
\eqtn
\elmm

Let $e'$ and $\tau'$ be as in Lemma~\ref{221020a}. Then by Lemma~\ref{170920c}, there is a bijection 
$$
f_{\tau'}:\Omega'\to S'_{\tau'},\ \alpha' \mapsto (x_{\alpha'},y_{\alpha'}),
$$
where $S'_{\tau'}$ and  $x_{\alpha'}$, $y_{\alpha'}$ are defined, respectively, by Eqs.~\eqref{200920w} and~\eqref{200920w3} for the scheme~$\cX'$. In what follows, we set $\varphi(r,s)=(r',s')$ for all $(r,s)\in S\times S$.

\lmml{151020a}
In the above notation, the mapping 
\qtnl{201020a}
f:=f_\varphi:\Omega\to\Omega',\ \alpha\mapsto [\varphi(\alpha^{f_\tau})]^{ f_{\tau'}^{-1}},
\eqtn
is a bijection. Moreover, $\tau^f=\tau'$ and given $\alpha,\beta\in\Omega$,
\qtnl{181020a}
\varphi(r(\alpha,\beta))=r(\alpha^f,\beta^f)
\eqtn 
whenever at least  one of the following conditions is satisfied:
\nmrt
\tm{a} $\alpha=\mu$,
\tm{b} $\alpha=\nu$  and $r(\mu,\alpha)\subseteq e$,
\tm{c} $\alpha=\rho$ and $r(\mu,\alpha)\nsubseteq e$.
\enmrt
\elmm
\proof The mappings $f_\tau:\Omega\to S_\tau$ and $f_{\tau'}^{-1}:S'_{\tau'}\to\Omega'$ are bijections by Lemma~\ref{170920c}. Since $\varphi$ induces a bijection $S^{}_{\tau^{}}\to S'_{\tau'}$, the mapping $f$ being the composition of these three bijections  is a bijection too.\medskip

Setting $r=r(\mu,\rho)$, by the first part of Eq.~\eqref{130920a} we have
\qtnl{210920a}
\mu^f=[\varphi(\mu^{f_\tau})]^{ f_{\tau'}^{-1}}=
[\varphi(1_{\Omega},r^*)]^{ f_{\tau'}^{-1}}=(1_{\Omega'},{r\phmb{*}}')^{ f_{\tau'}^{-1}}=
(1_{\Omega'},{r'}\phmb{*})^{ f_{\tau'}^{-1}}=\mu'.
\eqtn
Furthermore, setting  $\beta'=\beta^f$ and $(\beta')^{f_{\tau'}}=(x'_{\beta'},y'_{\beta'})$, we have
\qtnl{211020b}
(x'_{\beta'},y'_{\beta'})=(\beta')^{f_{\tau'}}=(\beta^f)^{f_{\tau'}}=
([\varphi(\beta^{f_\tau})]^{ f_{\tau'}^{-1}})^{f_{\tau'}}=
\varphi(\beta^{f_\tau})=((x^{}_{\beta^{}})',(y_{\beta^{}})').
\eqtn
It follows that $x'_{\beta'}=(x^{}_{\beta^{}})'$. By formula~(9),
$$
\varphi(r(\mu,\beta))=(x^{}_{\beta^{}})'=x'_{\beta'}=r(\mu',\beta')=r(\mu^f,\beta^f),
$$
which proves Eq.~\eqref{181020a} if the condition~(a) is satisfied. The same arguments work for  the conditions~(b) and~(c): the only difference is that one should use, respectively, the second and third parts of Eq.~\eqref{130920a}, and equality $y'_{\beta'}=(y^{}_{\beta^{}})'$, which follows from Eq.~\eqref{211020b}, instead of equality $x'_{\beta'}=(x^{}_{\beta^{}})'$.\eprf

\section{Proof of Theorem~\ref{4condition}}

In what follows, $\cX=(\Omega,S)$ is an  imprimitive pseudofrobenius scheme, $e$ is a nontrivial parabolic of~$\cX$, and $\fT(\cX)$ is defined as in Section~3.4 for $t=4$.

\lmml{140920a}
Assume that $\cX$ satisfies the $4$-condition, the arrays $T,T'\in\fT(\cX)$ are such that $T_{12}\subseteq e$, $T_{13} \nsubseteq e$, and 
\qtnl{270920a}
T^{}_{ij}=T'_{ij}~\text{for all}~\ 1\le i,j\le 4\ ~\text{except for}~\{i,j\}=
\css
\{2,4\} &\text{if $T_{14}\subseteq e$,}\\
\{3,4\} &\text{if $T_{14}\nsubseteq e$.}\\
\ecss
\eqtn
Then $T=T'$.
\elmm
\proof  We have $T=T(\gamma)$ and $T'=T(\gamma')$ for some $\gamma,\gamma'\in\Omega^4$. First, we prove that 
\qtnl{270920d}
(\gamma^{}_1,\gamma^{}_2)=(\gamma'_1,\gamma'_2)\quad\Rightarrow\quad \gamma=\gamma'.
\eqtn
Indeed condition~\eqref{270920a} implies that $T^{}_{13}=T'_{13}$ and $T^{}_{23}=T'_{23}$. Now if $(\gamma^{}_1,\gamma^{}_2)=(\gamma'_1,\gamma'_2)$ then $\gamma_3^{}$ and $\gamma_3'$ lie in  $\gamma_1 T_{13}\cap \gamma_2T_{23}$. On the other hand, $T_{12}\subseteq e$ and $T_{13} \nsubseteq e$ by the hypothesis. Therefore, $T^{}_{32}=r(\gamma_3,\gamma_2)=r(\gamma_1,\gamma_3)^*r(\gamma_1,\gamma_2)\subseteq T^*_{13}T^{}_{12}\nsubseteq e$ and then
$$
|\gamma_1 T_{13}\cap \gamma_2T_{23}|=c_{T_{13}T_{32}}^{T_{12}}=1
$$
by Lemma~\ref{intersection}. Thus, $\gamma^{}_3=\gamma'_3$. 

Further, if $T_{14}\not\subseteq e$, then the above argument with $\gamma^{}_3,\gamma'_3$ replaced by $\gamma^{}_4,\gamma'_4$  and $T_{13}$ and $T_{23}$ replaced by $T_{14}$ and $T_{24}$ shows that $\gamma_4^{}=\gamma'_4$. Now let $T_{14}\subseteq e$. By condition~\eqref{270920a}, we have $T^{}_{14}=T'_{14}$ and $T^{}_{34}=T'_{34}$. Since $\gamma_3^{}=\gamma_3'$, both $\gamma_4^{}$ and $\gamma_4'$ lie in  $\gamma_1 T_{14}\cap \gamma_3T_{34}$. By Lemma~\ref{intersection},
$$
|\gamma_1 T_{14}\cap \gamma_3T_{34}|=c_{T_{14}T_{43}}^{T_{13}}=1
$$
and hence again $\gamma_4^{}=\gamma'_4$. This completes the proof of Eq.~\eqref{270920d}.

Eq.~\eqref{270920d} shows that the set of all $\wt{\gamma}\in \Omega^4$ for which $(\wt{\gamma_1},\wt{\gamma_2})=(\gamma_1,\gamma_2)$ and $T(\wt{\gamma})=T$, is a singleton, namely $\{\gamma\}$. Since the scheme $\cX$ satisfies the $4$-condition, $c_T^r=1$ with $r=r(\gamma_1,\gamma_2)$. Taking into account that 
$$r(\gamma'_1,\gamma'_2)=T'_{12}=T^{}_{12}=r(\gamma_1,\gamma_2)=r,$$ 
we see that there exists a quadruple $\gamma''\in\fT(\cX)$ such that $T(\gamma'')=T$ and $(\gamma''_1,\gamma_2'')=(\gamma'_1,\gamma'_2)$. By Eq.~\eqref{270920d} for $\gamma=\gamma''$, we have $\gamma''=\gamma'$. Thus,
$$
T=T(\gamma'')=T(\gamma')=T',
$$
as required.\eprf\medskip

\lmml{120920a}
Let $\tau=(\mu,\nu,\rho)$ and $\tau'=(\mu',\nu',\rho')$ be base triples of~$\cX$, satisfying condition~\eqref{181020b} with respect to the algebraic isomorphism $\varphi=\id_S$ and bijection $f=f_\varphi$ defined by Eq.~\eqref{201020a}. Assume that the scheme $\cX$ satisfies the $4$-condition. Then $f\in\aut(\cX)$ and $\tau^f=\tau'$.
\elmm
\proof To prove that $f\in\aut(\cX)$, it suffices to verify that Eq.~\eqref{181020a} holds for all distinct $\alpha,\beta\in\Omega$. Note that 
if $\{\alpha,\beta\}\subseteq \{\mu,\nu,\rho\}$, then this follows from  Lemma~\ref{151020a}.\medskip
 
{\bf Claim.} {\it Eq.~\eqref{181020a} holds if $\alpha$ or $\beta$ belongs to~$\{\mu,\nu,\rho\}$.} \medskip

\proof Without loss of generality, we assume that $\alpha\in\{\mu,\nu,\rho\}$.  Let $T=T(\gamma)$ and  $T'=T(\gamma')$, where $\gamma=(\mu,\nu,\rho,\beta)$ and $\gamma'=(\mu',\nu',\rho',\beta^f)$. Then by above, 
$$
T^{}_{ij}=T'_{ij},\quad 1\le i,j\le 3,
$$
and also
$$
T^{}_{12}\subseteq e\qaq T^{}_{13}\nsubseteq e,
$$
because $\tau$ and $\tau'$ are base triples.
Finally, by Lemma~\ref{151020a}, we have
$$
T^{}_{14}=T'_{14}\qaq
\css
T^{}_{34}=T'_{34} & \text{if $T_{14}\subseteq e$},\\
T^{}_{24}=T'_{24} & \text{if $T_{14}\nsubseteq e$}.\\
\ecss
$$
Thus all conditions of Lemma~\ref{140920a} are satisfied and hence $T=T'$. Since $\alpha=\gamma^{}_i$ and $\alpha^f=\gamma'_i$ for some $1\le i\le 3$, we conclude that
$$
r(\alpha^f,\beta^f)=T'_{i4}=T^{}_{i4}=r(\alpha,\beta),
$$
as required.\eprf\medskip

To complete the proof of Lemma~\ref{120920a}, let $T=T(\gamma)$, where this time $\gamma=(\mu,\alpha, \nu, \beta)$ if $(\mu,\beta)\in e$, and $\gamma=(\mu,\nu,\alpha,\beta)$ if $(\mu,\beta)\notin e$. Then by Claim, all conditions of Lemma~\ref{140920a} are satisfied for $T'=T(\gamma')$ with $\gamma'=\gamma^f$. Thus, $T=T'$. It follows  that
$$
r(\alpha^f,\beta^f)=T'_{ij}=T^{}_{ij}=r(\alpha,\beta),
$$
where $(i,j)=(2,4)$ if $(\mu,\beta)\in e$ and $(i,j)=(3,4)$ if $(\mu,\beta)\notin e$.
\eprf\medskip

{\bf Proof of Theorem~\ref{4condition}.} The necessity follows from Section~3.4. To prove the sufficiency, let~$\cX$ be an imprimitive pseudofrobenius scheme, and $e$ a nontrivial parabolic of~$\cX$. Assume that $\cX$ satisfies the $4$-condition. In view of~\cite[Corollary~3.3.9]{CP}, every schurian imprimitive equivalenced scheme is Frobenius. Thus, it suffices to prove that $\mathcal{X}$ is schurian, i.e., the group $\aut(\cX)$ acts transitively on each~$s\in S$.\medskip

Assume first that $s\subseteq e$, and $(\mu,\nu), (\mu',\nu')\in s$. Take arbitrary point $\rho$ so that $\tau=(\mu,\nu,\rho)$ is a base triple of $\cX$ with respect to~$e$. Setting $r=r(\mu,\rho)$ and $t=r(\nu,\rho)$, we have $c_{st}^r\ne 0$. It follows that there exists a point~$\rho'$ such that $\tau'=(\mu',\nu',\rho')$ is a base triple of $\cX$ with respect to~$e$, such that condition~\eqref{181020b} is satisfied. By Lemma~\ref{120920a}, we have $(\mu,\nu)^f=(\mu',\nu')$, where $f$ is the  automorphism of~$\cX$, defined in that lemma. Thus, $\aut(\cX)$ acts transitively on~$s$.\medskip

Now let $s\nsubseteq e$ and $(\mu,\rho),(\mu',\rho')\in s$. Take arbitrary point $\nu$ so that $\tau=(\mu,\nu,\rho)$ is a base triple of $\cX$ with respect to~$e$. Then as above, there exists a point~$\nu'$ such that $\tau'=(\mu',\nu',\rho')$ is a base triple of $\cX$ with respect to~$e$, such that condition~\eqref{181020b} is satisfied. Using Lemma~\ref{120920a} again, we find $f\in\aut(\cX)$ such that $(\mu,\nu)^f=(\mu',\nu')$. Thus, $\aut(\cX)$ acts transitively on~$s$.\eprf


\section{Proof of Theorem~\ref{2}}

A separable pseudofrobenius scheme is obviously Frobenius. So every proper pseudofrobenius scheme is not separable. Therefore, Theorem~\ref{2} follows from the theorem below.

\thrml{separability}
An imprimitive pseudofrobenius scheme $\cX$ is separable unless the conclusion of Theorem~\ref{2} holds for some (and hence for all) Frobenius group~$G$ such that $\cX$ is $G$-pseudofrobenius.
\ethrm

Throughout the rest of the section, $\mathcal{X}=(\Omega,S)$ and~$n$ and~$k$ the degree and valency of $\cX$, respectively. We need an auxiliary lemma. 

\lmml{3parabolics}
Suppose that $\mathcal{X}$ has a chain of parabolics $1_{\Omega}<e_1<e_2\leq e_3<\Omega\times \Omega$. Then $\mathcal{X}$ is separable unless $e_2=e_3$ and
\qtnl{multisets} 
\{\{n_1-1,\frac{n_2}{n_1}-1,\frac{n}{n_2}-1\}\}=\{\{k,k,k\}\}~\text{or}~\{\{k,k,2k\}\},
\eqtn
where $n_1=n_{e_1}$ and $n_2=n_{e_2}$.
\elmm

\proof
By~\cite[Theorem~5.1]{P}, every scheme of degree~$n$ and maximal valency~$k$ is separable whenever $n>3ck(k-1)$, where $c$ is the indistinguishing number of this scheme. By definition, $c$ is the sum of some intersection numbers, and $c=k-1$ if the scheme is Frobenius~\cite[Lemma~2.2]{MP} and hence if it is pseudofrobenius. Thus, the scheme $\mathcal{X}$ is separable if 
\qtnl{sufficient}
n>3k(k-1)^2.
\eqtn

By Lemma~\ref{divide}, the number $k$ divides $n_1-1$, $\frac{n_2}{n_1}-1$, $\frac{n_3}{n_2}-1$, and $\frac{n}{n_3}-1$. If $e_2\neq e_3$ then each of these numbers is nonzero and hence greater than or equal to~$k$. Consequently,
$$n=n_1\frac{n_2}{n_1}\frac{n_3}{n_2}\frac{n}{n_3}>(n_1-1)\left(\frac{n_2}{n_1}-1\right)\left(\frac{n_3}{n_2}-1\right)\left(\frac{n}{n_3}-1\right)\geq k^4>3k(k-1)^2.$$ 
On the other hand, let $e_2=e_3$. Suppose that Eq.~\eqref{multisets} does not hold. Then the above divisibility argument easily shows that 
$$n=n_1\frac{n_2}{n_1}\frac{n}{n_2}>(n_1-1)\left(\frac{n_2}{n_1}-1\right)\left(\frac{n}{n_2}-1\right)\geq 3k^3>3k(k-1)^2.$$
Thus, in both cases Eq.~\eqref{sufficient} holds in both cases and $\mathcal{X}$ is separable. 
\eprf\medskip

{\bf Proof of Theorem~\ref{separability}.}
Let $G$ be a Frobenius group such that~$\mathcal{X}$ is $G$-pseudofrobenius, and $H$ and $K$ the kernel and complement of $G$, respectively. Clearly, $|H|=n$ and $|K|=k$. Assume that $\mathcal{X}$ is not separable.

\lmml{3groups}
Suppose that $H$ has a chain of $K$-invariant subgroups $\{1\}<H_1<H_2\leq H_3<H$. Then $H_2=H_3$ and Eq.~\eqref{multisets} holds for $n_1=|H_1|$ and $n_2=|H_2|$.
\elmm

\proof
By Lemma~\ref{parabolics}(2), the condition of Lemma~\ref{3parabolics} holds for $e_i=e(H_i)$ and also $n_{e_i}=n_i$  for $i=1,2,3$. Thus, the required statement is a consequence of Lemma~\ref{3parabolics}.~
\eprf\medskip

Let $\pi=\pi(H)$ and $d=d(G)$.

\lmml{pr1}
$|\pi|\leq 3$ and $d\leq 3$.
\elmm

\proof
By Lemma~\ref{3groups}, it suffices to show that if $|\pi|\geq 4$ or $d\geq 4$, then $H$ has a chain of $K$-invariant subgroups $\{1\}<H_1<H_2<H_3<H$. This is obvious if $d\geq 4$. Now let $|\pi|\geq 4$ and let $P$, $Q$, $R$ be pairwise distinct Sylow subgroups of~$H$. Put $H_1=P$, $H_2=P\times Q$, and $H_3=P\times Q \times R$. Since $H_1$, $H_2$, and $H_3$ are $K$-invariant, we are done.
\eprf\medskip

\lmml{pr2}
$|\pi|\leq 2$.
\elmm

\proof
By Lemma~\ref{pr1}, we may assume towards to contradiction that $|\pi|=3$. Let $\syl(H)=\{P,Q,R\}$. Then the condition of Lemma~\ref{3groups} holds for $H_1=P$ and $H_2=H_3=P\times Q$. Moreover, 
$$\left|\left\{n_1-1,\frac{n_2}{n_1}-1,\frac{n}{n_2}-1\right\}\right|=\left|\left\{|P|-1,|Q|-1,|R|-1\right\}\right|=3,$$
which contradicts Lemma~\ref{3groups}.
\eprf\medskip

Note that $d\geq 2$, since the scheme $\mathcal{X}$ is imprimitive. Together with Lemmas~\ref{pr1} and~\ref{pr2}, this implies that $(|\pi|,d)\in \{1,2\}\times \{2,3\}$. Further, let $d=3$ and 
$$\{1\}<H_1<H_2<H$$ 
a chain of $K$-invariant subgroups. The condition of Lemma~\ref{3groups} is satisfied for $H_1$ and $H_2=H_3$ and Eq.~\eqref{multisets} holds. To complete the proof, we consider two cases depending on $|\pi|=1$ or~$2$. In what follows,~$S$ is a factor block such that the section~$G^S$ is principal; without loss of generality, we may assume that $S=H_1$, $H_2/H_1$, or $H/H_2$.\medskip

{\bf Case~1:} $|\pi|=1$. In this case, $|H|$ is a $p$-group for a prime~$p$. In particular, each of the numbers $n_1$, $\frac{n_2}{n_1}$, $\frac{n}{n_2}$ is a $p$-power. It follows that the ratio of every two elements from the set $\{n_1-1,\frac{n_2}{n_1}-1,\frac{n}{n_2}-1\}$ is other than~$2$. Therefore, 
\qtnl{set1} 
\left\{n_1-1,\frac{n_2}{n_1}-1,\frac{n}{n_2}-1\right\}=\{k\}
\eqtn
by Eq.~\eqref{multisets} and hence 
$$|H|=n=n_1\frac{n_2}{n_1}\frac{n}{n_2}=(k+1)^3;$$
in particular, $|H|=p^{3a}$ for an integer~$a\geq 1$, and $|S|=|H_1|=|H_2/H_1|=|H/H_2|=k+1=p^a$. Due to Eq.~\eqref{set1}, the action of $K\leq G$ on $S$ is transitive on $S\setminus \{1\}$. Thus, $G^S$ is $2$-transitive.\medskip

{\bf Case~2:} $|\pi|=2$. In this case, $\pi=\{p,q\}$ for distinct primes $p$ and $q$. By Lemma~\ref{class}(3), the set $\{n_1,\frac{n_2}{n_1},\frac{n}{n_2}\}$ contains both $p$- and $q$-powers and hence has cardinality at least~$2$. Together with Eq.~\eqref{multisets}, this implies that 
\qtnl{set2} 
\left\{\left\{n_1-1,\frac{n_2}{n_1}-1,\frac{n}{n_2}-1\right\}\right\}=\{\{k,k,2k\}\}.
\eqtn
Without loss of generality, we may assume that $k=p^a-1$ and $2k=q^b-1$ for some integers $a,b\geq 1$; in particular, $q^b=2p^a-1$. By Eq.~\eqref{set2},  
$$|H|=n=n_1\frac{n_2}{n_1}\frac{n}{n_2}=(k+1)^2(2k+1)=p^{2a}q^b,$$
and $|S|\in \{n_1-1,\frac{n_2}{n_1}-1,\frac{n}{n_2}-1\}=\{p^a,q^b\}$. Then the action of $K\leq G$ on $S$ has one (if $|S|=p^a$) or two (if $|S|=q^b$) orbits on $S\setminus \{1\}$. In the first case $G^S$ is $2$-transitive whereas in the second one it is of rank~$3$.
\eprf\medskip

{\bf Proof of Corollary~\ref{nilpotent}.}
We keep notation from Theorem~\ref{2}. By this theorem $d\leq 3$ and hence $c\leq d\leq 3$ by Lemma~\ref{class}(1).\medskip

{\bf Claim~1.} If $d=3$ and $2\notin \pi$, then $H$ is abelian.\medskip

\proof
Since $d=3$, the group $G$ has a $2$-transitive principal section of degree $p^a$ by Theorem~\ref{2}. This implies that $|K|=k=p^a-1$. Then $|K|$ is even because $p\neq 2$. Therefore, $H$ is abelian by Lemma~\ref{class}(4).
\eprf\medskip

{\bf Claim~2.} If $|\pi|=2$, then $c\leq d-1$.\medskip

\proof
Let $\syl(H)=\{P,Q\}$, where $P\neq Q$. In this case $H=P\times Q$ and
$$c\leq \max\{c(P),c(Q)\}\leq \max\{d(G^P),d(G^Q)\}\leq d-1,$$
where the second and third inequalities hold statements~(1) and~(2) of Lemma~\ref{class}, respectively.
\eprf\medskip

Now suppose that $H$ is not a $2$-group or $d\neq 3$. If $(|\pi|,d)=(1,3)$, then $H$ is abelian by Claim~$1$; if $(|\pi|,d)=(1,2)$, then $c\leq 2$ by Lemma~\ref{class}(1); if $(|\pi|,d)=(2,2)$, then $H$ is abelian by Claim~$2$; if $(|\pi|,d)=(2,3)$ and $2\notin \pi$, then $H$ is abelian by Claim~$2$; if $(|\pi|,d)=(2,3)$ and $2\in \pi$, then $c\leq 2$ by Claim~$2$. Thus, statement~(1) of the corollary holds, whereas statement~(2) immediately follows from statement~(1).
\eprf

\section{Proof of Theorem~\ref{3}}

It suffices to prove the ``only if'' part. Assume that the groups $G_1$ and $G_2$ have the same character tables. Then there are bijections between the irreducible characters and the conjugacy classes of $G_1$ and $G_2$, which preserve the degrees and cardinalities, respectively. In particular, this defines a one-to-one correspondence between the basis relations of the (commutative) schemes~$\mathcal{Y}_1$ and~$\mathcal{Y}_2$ of the groups $G_1\Inn(G_1)\leq \sym(G_1)$ and $G_2\Inn(G_2)\leq \sym(G_2)$, which preserves the intersection numbers (see, e.g.,~\cite[Chapter~II, Theorems~3.6 and~7.2]{BI}). Thus, $\mathcal{Y}_1$ and $\mathcal{Y}_2$ are algebraically isomorphic.\medskip

\lmml{camina}
The schemes $\mathcal{X}_1=\inv(G_1)$ and $\mathcal{X}_2=\inv(G_2)$ are algebraically isomorphic.
\elmm

\proof
By~\cite[Theorem~3.4.8(1)]{CP}, it suffices to prove that the scheme $\mathcal{Y}=\mathcal{Y}_i$ is isomorphic to the wreath product of the scheme $\mathcal{X}=\mathcal{X}_i$ and some scheme, $i=1,2$. Note that $\mathcal{X}$ and $\mathcal{Y}$ are Cayley schemes over the kernel $H$ of $G$ and $G$, respectively. In view of~\cite[Theorem~2.4.16]{CP}, we will prove the required statement in the language of Schur rings.\medskip

Let $\mathcal{A}$ be the Schur ring over $G$ defined by the group $G\Inn(G)\leq \sym(G)$; the basic sets of~$\mathcal{A}$ are just the conjugacy classes of $G$. The group $H$ being a normal abelian subgroup of $G$ is a union of some basic sets of $\mathcal{A}$. Moreover, the restriction~$\mathcal{A}_{H}$ of~$\mathcal{A}$ on~$H$ is the Schur ring defined by $G\leq \sym(H)$. On the other hand, each conjugacy class of $G$ outside $H$ is a union of some $H$-cosets by~\cite[p.~153]{Camina}. This means that $XH=HX=X$ for each basic set $X\subseteq G\setminus H$. Thus, $\mathcal{A}$ is the wreath product of $\mathcal{A}_H$ and some Schur ring over $G/H$. 
\eprf\medskip

By the hypothesis of the theorem, the conclusion of Theorem~\ref{2} does not hold for $G_1$. It follows that $\mathcal{X}_1$ is separable (Theorem~\ref{separability}). Therefore, $\mathcal{X}_1$ and $\mathcal{X}_2$ are isomorphic (Lemma~\ref{camina}). In particular, the groups $\aut(\mathcal{X}_1)$ and $\aut(\mathcal{X}_2)$ are isomorphic. Since $G_1$ and $G_2$ are imprimitive Frobenius groups, $G_1=\aut(\mathcal{X}_1)$ and $G_2=\aut(\mathcal{X}_2)$ by~\cite[Theorem~2.5.8]{FKM}. Thus, $G_1$ and $G_2$ are isomorphic.

\section{Proof of Theorem~\ref{circulant}}
Denote by~$H$ the cyclic group underlying the circulant~$\fX$. Then 
$$H\leq \aut(\fX)=:G$$ 
(we identify $H$ with the image of its regular representation). Since $G$ is a Frobenius group and every nontrivial element of $H$ is fixed point free, $H$ is the kernel of $G$. In particular, $H$ is normal in $G$. 

According to~\cite[Section~2.6.1]{CP}, there exists the minimal scheme $\cX$ such that $\aut(\cX)=\aut(\fX)$. In particular, $\cX$ is a normal Cayley scheme (of~$H$) in the sense of~\cite[Section~2.4]{CP}. Moreover, in view of~\cite[Exercise~4.7.37(1)]{CP}, we have
$$\cX=\inv(\aut(\cX))=\inv(\aut(\fX))=\inv(G).$$
Consequently, $\cX$ is a Frobenius scheme.

Let $n\notin \{p,p^2,p^3,pq,p^2q\}$, where $p$ and $q$ are distinct primes. Then 
$$|\pi(n)|\geq 3\quad\text{or}\quad\Omega(n)\geq 4,$$
where $\pi(n)$ is the set of all prime divisors of~$n$ and $\Omega(n)$ is the total number of prime divisors of~$n$ (with multiplicities). Clearly, $|H|=n$ and $\pi(H)=\pi(n)$. Since $H$ is cyclic, every subgroup of $H$ is characteristic and the length of the longest chain of subgroups of~$H$ equals $\Omega(n)$. Therefore, 
$$|\pi(H)|\geq 3\quad\text{or}\quad d(G)\geq 4.$$
In both cases, the scheme $\cX$ is separable by Theorem~\ref{separability}. Together with~\cite[Theorem~2.5]{FKV}, this implies that $\dimwl(\fX)\leq 2$. The classification of all regular graphs of the Weisfeiler-Leman dimension~$1$~\cite[Lemma~3.1 (a)]{AKRV} yields~$\dimwl(\fX)\neq 1$. Thus, $\dimwl(\fX)=2$.  
\eprf


\begin{thebibliography}{1}

\bibitem{Andre}
J.~Andr\'{e}, \emph{Eine geometrische Kennzeiehnung imprlmitiver Frobeniusgruppen}, Abh. Math. Sem. Univ. Hamburg, {\bf 51}, 120--135 (1981).

\bibitem{AKRV}
V.~Arvind, J.~K\"obler, G.~Rattan, and O.~Verbitsky, \emph{Graph Isomorphism, Color Refinement, and Compactness}, Comput. Complexity, \textbf{26}, No.~3, 627--685 (2017).

\bibitem{BI}
E.~Bannai and T.~Ito, \emph{Algebraic Combinatorics I. Association schemes}, The Benjamin/Cummings Publishing Co., Inc., Menlo Park, CA (1984).

\bibitem{BCN}
A.~Brouwer, A.~Cohen, and A.~Neumaier, \emph{Distance-regular graphs}, Springer, Heidelberg, (1989).

\bibitem{Bourbaki}
N.~Bourbaki, \emph{Elements of mathematics. Algebra I. Chapters~1-3}, Springer-Verlag (1998).


\bibitem{Camina}
A.~R.~Camina, \emph{Some conditions which almost characterize Frobenius groups}, Isr. J. Math., {\bf 31}, No.~2, 153--160 (1978).

\bibitem{CP}	
G.~Chen and I.~Ponomarenko, \emph{Coherent Configurations}, Central China Normal University Press (2019); the updated version is available at \url{http://www.pdmi.ras.ru/~inp/ccNOTES.pdf}.



\bibitem{CHPV}
G.~Chen, J.~He, I.~Ponomarenko, and A.~Vasil'ev, \emph{A characterization of exceptional pseudocyclic association schemes by multidimensional intersection numbers}, Ars Math. Contemp., {\bf 21}, No.~1, \#P1.10 (2021).


\bibitem{EP}
S.~Evdokimov and I.~Ponomarenko, \emph{On primitive cellular algebras}, J. Math. Sci. (N.-Y.), {\bf 107}, No.~5, 4172--4191 (2001). 


\bibitem{EP2}
S.~Evdokimov and I.~Ponomarenko, \emph{Schurity of $S$-rings over a cyclic group and generalized wreath product of permutation groups}, St. Petersburg Math. J., {\bf 24}, No.~3, 431--460 (2013).

\bibitem{FKM}
I.~Farad\u{z}ev, M.~Klin, and M.~Muzychuk, \emph{Cellular rings and groups of automorphisms of graphs}, In: I.~Farad\u{z}ev et al (eds.), Investigations in Algebraic Theory of Combinatorial Objects, Springer Science+Business Media, 1--152 (1994).


\bibitem{FKV}
F.~Fuhlbr\"uck, J~ K\"obler,  and O.~Verbitsky, \emph{Identiability of graphs with small color classes by the Weisfeiler-Leman algorithm}, in: Proc. $37$th International Symposium on Theoretical Aspects of Computer Science, Dagst\"uhl Publishing, Germany, 43:1--43:18 (2020).


\bibitem{GN}
M.~Gr\"ohe and  D.~Neuen, \emph{Recent Advances on the Graph Isomorphism Problem}, http://arxiv.org/abs/2011.01366 [cs.DS], 1--39 (2020).


\bibitem{HH}
 D.~G.~Higman, \emph{Characterization of families of rank $3$ permutation groups by the subdegrees. II}, Arch. Math. {\bf  21}, 353--361 (1970).


\bibitem{Is}
I.~M.~Isaacs, \emph{Finite group theory}, Graduate studies in Mathematics, {\bf 92}, American Mathematical Society (2008).







\bibitem{MP} 
M.~Muzychuk and  I.~Ponomarenko, \emph{On pseudocyclic association schemes}, Ars Math. Contemp., {\bf 5}, No.~1, 1--25 (2012).





\bibitem{P} 
I.~Ponomarenko, \emph{On the separability of cyclotomic schemes over finite field}, Algebra Analiz, {\bf 32}, No.~6, 124--146 (2020).



\bibitem{Ry}
G.~Ryabov, \emph{On separable abelian $p$-groups}, Ars Math. Contemp., {\bf 17}, No.~2, 467--479 (2019).

\bibitem{Suz}
M.~Suzuki, \emph{A new type of simple groups of finite order}, Proc. Natl. Acad. Sci. U.S.A., {\bf 46}, No.~6, 868--870 (1960).
 

\bibitem{Wi}
H.~Wielandt, \emph{Finite permutation groups}, Academic Press, New York - London (1964).

	
\end{thebibliography}
\end{document}